\documentclass{article}
\usepackage{amsfonts}

\begin{document}

\def\Tr{\mbox{Tr \-}}
\def\inj{\mbox{inj \-}}
\def\R{{\mathbb R}}
\def\C{{\mathbb C}}
\def\H{{\mathbb H}}
\def\Ca{{\mathbb Ca}}
\def\Z{{\mathbb Z}}
\def\N{{\mathbb N}}
\def\Q{{\mathbb Q}}
\def\Ad{\mbox{Ad \-}}
\def\k{{\bf k}}
\def\l{{\bf l}}
\def\sp{\mbox{\bf sp}}

\title{Special K\"ahler Metrics on Complex Line Bundles
             and the Geometry of $K3$-Surfaces}
\author{Ya. V. Bazaikin}
\date{}
\maketitle

%\thanks      The author was supported by
%            the Russian Foundation for Basic Research
%            (Grant~03--01--00403),
%            the RAS Program
%             ``Mathematical Methods in Nonlinear Dynamics,''
%            and
%            the Program ``Development of the Scientific Potential of Higher School''
%            of the Ministry for Education of the Russian Federation
%            (Grant~8311).

%\keywords    special K\"ahler manifold, $K3$-surface
%\endkeywords
%on the tangent bundles of weighted complex projective lines. We give
%a~geometric description of a~neighborhood of the moduli space of
%special K\"ahler metrics on a~$K3$-surface.
%\endabstract
%\endtopmatter

\section[]{Introduction}

In this article we continue studying
the Ricci-flat Riemannian metrics that were constructed in~[1].
On closer examination it turned out that they possess
a~number of remarkable properties; in particular,
they have the holonomy group~$SU(2)$,
so presenting special K\"ahler metrics.

The metrics of holonomy~$SU(2)$ are interesting because of
their applications
in mathematical physics.
In superstring theory and $M$-theory
there appear compact manifolds
with special holonomy groups.
Moreover, if we admit the presence of physically isolated
singularities then it suffices to study asymptotically
flat metrics on the normal bundles of these singularities.
Thus, we arrive at the problem of studying asymptotically locally
Euclidean metrics with special holonomies on bundles over orbifolds.

One of the most topical examples of special K\"ahler metrics is the
Eguchi--Hanson metric~[2] on the cotangent bundle $T^*S^2$ of the
standard two-dimensional sphere (without singularities). The
Eguchi--Hanson metric has played an~important role in studying
special holonomy groups. Namely, Page~[3] proposed a~description of
the space of special K\"ahler metrics on a~$K3$-surface in which the
Eguchi--Hanson metric plays the role of an~``elementary brick.''
More exactly, represent a~$K3$-surface using Kummer's construction;
i.e., consider the involution of the flat torus $T^4$ which arises
from the central symmetry of the Euclidean space~${\Bbb R}^4$.
Factorizing, we obtain an~orbifold with $16$ singular points whose
neighborhoods look like~${\Bbb C}^2/{\Bbb Z}_2$. Blowing up the
resulting orbifold in a~neighborhood of each singular point, we
obtain a~$K3$-surface. Topologically, the construction of blowing up
a~singular point of the form ${\Bbb C}^2/{\Bbb Z}_2$ is carried out
as follows: We have to delete the singularity and identify its
neighborhood with the space of the spherical bundle in~$T^*S^2$
without the zero fiber~$S^2$. Page proposed to consider a~metric
on~$T^*S^2$ which is homothetic to the Eguchi--Hanson metric with
a~sufficiently small homothety coefficient so that the metric on the
boundary of the glued spherical bundle becomes arbitrarily close to
a~flat metric. After that we need to deform slightly the metric on
the torus so as to obtain a~smooth metric on a~$K3$-surface with
holonomy~$SU(2)$. A~simple evaluation of the degrees of freedom in
the process of this operation demonstrates that we obtain
a~$58$-dimensional family of metrics which agrees with the
well-known results on the dimension of the moduli space of such
metrics~[4]. Later, Page's idea was used by Joyce~[5,\,6] for
constructing the first compact examples of manifolds with the exotic
holonomy groups~$G_2$ and~$\mbox{Spin}(7)$.

We can identify the two-dimensional sphere $S^2$ with
the complex projective line~${\Bbb C} P^1$.
Consider its natural generalization ${\Bbb C} P^1(k,l)$,
the weighted complex projective line,
which is a~complex orbifold with two singular points.
In this article we obtain the following

\vskip0.2cm

{\bf Theorem.} {\it On the cotangent bundle $M_{k,l}=T^* {\Bbb C}
P^1(k,l)$ of the weighted complex projective line, there is a~metric
with the holonomy group~$SU(2)$.}

\vskip0.2cm

This metric was found in~[1] in a~special coordinate system
as a~solution to the equation of the zero Ricci curvature on
torus bundles over two-dimensional surfaces.
Also, in~[1] it was proven that the metric has the isometry group
$U(1) \times U(1)$ and cohomogeneity ~2 except for the case $k=l=1$.
For $k=l=1$ our metric coincides with the Eguchi--Hanson metric.

Asymptotically, the constructed metrics behave as follows:
At infinity the metric tends to the Euclidean metric on~${\Bbb C}^2/{\Bbb Z}_{k+l}$,
while in a~neighborhood of each of the two singular points
it tends to the respective Euclidean metrics on~
${\Bbb C}^2/{\Bbb Z}_k$ and ${\Bbb C}^2/{\Bbb Z}_l$.
Therefore, with Page's idea in mind, we propose to use
the metrics on~$M_{k,l}$
for blowing up
the singularities of the
form ${\Bbb C}^2/{\Bbb Z}_p$ on the orbifolds with holonomy~$SU(2)$
in several steps: successively replace
each singularity with two singularities of less order
gluing the space $M_{k,l}$ with the constructed metric;
hopefully, we eventually ``remove'' all singular points.

As an~application, we consider a~representation of a~$K3$-surface
as the blow-up
of the singularities of the orbifold~$T^4/{\Bbb Z}_p$
for a~prime $p \neq 2$.
It turns out that the only possible case is $p=3$
in which we have to blow up $9$ singular points of the form
${\Bbb C}^2/{\Bbb Z}_3$.
This is done in two steps: first, using $M_{1,2}$, we
obtain $9$ singular points of the form
${\Bbb C}^2/{\Bbb Z}_2$ and then remove them by means of
$M_{1,1}=T^*S^2$.
Each time we slightly deform the metric on the ``glued'' space
and eventually obtain a~$K3$-surface with a~family of metrics
with holonomy~$SU(2)$.
Simple calculations demonstrate that the dimension of the so-obtained family
equals~$58$, as expected.
However, the so-described metrics on a~$K3$-surface differ essentially
from the family constructed by Page: actually, we give an~asymptotic
description of the moduli space of the metrics with holonomy
$SU(2)$ in a~neighborhood of a~flat metric on~$T^4/{\Bbb Z}_3$,
while Page gave an~asymptotic description of the same
space in a~neighborhood of a~flat metric on~$T^4/{\Bbb Z}_2$.
To justify  our construction rigorously, we establish
a~connection between the constructed metrics and multi-instantons~[7,\,8];
our metrics are the limit case of a~multi-instanton
corresponding to two sources with different ``masses.''

In the next section we consider weighted complex projective spaces
and describe the structure of~$M_{k,l}$ and the metric on it.
In~\S\,3 and \S\,4 we discuss applications to the geometry of $K3$-surfaces.

{\bf Acknowledgement}. The author is grateful to I.~A. Ta\u\i manov
for useful discussions. This work was supported by the Russian
Foundation for Basic Research (Grant 03-01-00403), the RAS Program
"Mathematical Methods in Nonlinear Dynamics," and the Program
"Development of the Scientific Potential of Higher School" of the
Ministry for Education of the Russian Federation (Grant 8311).

\section[]{A~Special K\"ahlerian Structure
on~$M_{k,l}$}

Let ${\bold k}=(k_0, k_1, \dots, k_n)$ be a~set of coprime positive integers.
Consider the action of ${\Bbb C}^*={\Bbb C} \backslash \{ 0 \}$ on
${\Bbb C}^{n+1} \backslash \{0\}$
with weights $k_0, \dots, k_n$:
$$
\lambda \in {\Bbb C}^* : (z_0, z_1, \dots, z_n) \mapsto (z_0
\lambda^{k_0}, z_1 \lambda^{k_1}, \dots, z_n \lambda^{k_n}).
$$
The orbit space of this action of~${\Bbb C} P^n({\bold k})$
possesses the structure of a~complex orbifold and is called
the~{\it weighted complex projective space}.
We denote the orbit of a~point $(z_0, \dots, z_n)$ by~$[z_0: \dots :z_n]$.
The structure of singularities of the orbifold
${\Bbb C} P^n({\bold k})$
can be rather complicated and depends in general on the
mutual divisibility
of different sets of the numbers~$k_0, \dots, k_n$.
In the case when each pair~$k_i, k_j$ is coprime
(this is the case we are interested in), the situation becomes somewhat simpler
and ${\Bbb C} P^n({\bold k})$ possesses only a~discrete collection of isolated singularities
$[1:0:\dots :0], [0:1:\dots :0], \dots, [0:0:\dots :1]$
beyond which it is a~complex analytic manifold.
As a~uniformizing atlas we have to consider the collection of~ charts:
$$
\phi_i (z_1, \dots, z_n) = [z_1: \dots : z_i : 1: z_{i+1} : \dots
: z_n], \quad  i=0,\dots , n.
$$
Moreover, for each chart, the uniformizing group is the group
$\Gamma_i={\Bbb Z}_{k_i}$
generated by the element
$\omega_i=e^{\frac{2 \pi i}{k_i}}$ whose action is given as follows:
$$
\omega_i (z_1, \dots, z_n) = \bigl(z_1 \omega_i^{k_0}, \dots , z_i
\omega_i^{k_{i-1}}, z_{i+1} \omega_i^{k_{i+1}}, \dots , z_n
\omega_i^{k_n} \bigr).
$$

We need a~generalization of the above construction.
Suppose that we have collections
${\bold k}=(k_0, \dots, k_n)$
and
${\bold l}=(l_0, \dots, l_m)$
of pairwise coprime integers such that $k_i>0$ and $l_j<0$.
As above, we can formally consider the action of~${\Bbb C}^*$ on
${\Bbb C}^{n+m+2}\backslash \{0\}$
with weights $k_0, \dots, k_n, l_0, \dots, l_m$.
Moreover, as the orbit space
${\Bbb C}P^{n+m+1}({\bold k},{\bold l}) =({\Bbb C}^{n+m+2}\backslash \{0\})/{\Bbb C}^*$
we obtain a~topological space possessing a~uniformizing atlas which is
given by the collection of charts $\phi_i$ and $\psi_j$,
$i=0,\dots, n$, $j=0, \dots, m$, defined as above.
However, the so-obtained orbit space does not possess the Hausdorff property:
the singular points corresponding to positive weights are not separated from
those corresponding to negative weights.
More exactly, consider the obvious embeddings of
${\Bbb C}^{n+1}={\Bbb C}^{n+1}\times\{0\}$ and
${\Bbb C}^{m+1}=\{0\}\times {\Bbb C}^{m+1}$ into
${\Bbb C}^{n+m+2}={\Bbb C}^{n+1}\times{\Bbb C}^{m+1}$.
Passing to the orbit space, we obtain two orbifolds
${\Bbb C} P^{n}({\bold k})$ and ${\Bbb C} P^m ({\bold l})$
that are naturally embedded into
${\Bbb C} P^{n+m+1}({\bold k},{\bold l})$
and cover together all singular points.
It is obvious that ${\Bbb C} P^n({\bold k})$ cannot be
separated
from ${\Bbb C} P^m({\bold l})$ in the quotient topology of the orbit space.
Therefore, define the two weighted complex projective spaces:
$$
{\Bbb C} P^{n+m+1}_+ ({\bold k},{\bold l})
= {\Bbb C} P^{n+m+1}({\bold k},{\bold l}) \backslash {\Bbb C} P^m ({\bold l}),
$$
$$
{\Bbb C} P^{n+m+1}_- ({\bold k},{\bold l})
= {\Bbb C} P^{n+m+1}({\bold k},{\bold l}) \backslash {\Bbb C} P^n ({\bold k}).
$$
It is clear that
${\Bbb C} P^{n+m+1}_+ ({\bold k},{\bold l})$ and
${\Bbb C} P^{n+m+1}_- ({\bold k},{\bold l})$ are noncompact
orbifolds with uniformizing atlases given by the respective collections of charts
$\phi_i$, $i=0,\dots,n$, and $\psi_j$, $j=0,\dots, m$.
Note that there is an~obvious isomorphism between the complex manifolds
${\Bbb C} P^{n+m+1}_+({\bold k},{\bold l}) \backslash {\Bbb C} P^n ({\bold k})$
and
${\Bbb C} P^{n+m+1}_-({\bold k},{\bold l}) \backslash {\Bbb C} P^m ({\bold l})$
induced by the identical transformation of the space~${\Bbb C}^{n+m+2}$.

We turn to description of some necessary spaces.
Consider a~pair $k$, $l$ of coprime positive numbers.
We obtain the weighted complex projective line
$S^2(k,l)={\Bbb C} P^1(k,l)$.
Let the uniformizing atlas
on~$S^2(k,l)$ consist of the two charts
$\phi_0$ and~$\phi_1$ defining the coordinates $z \in {\Bbb C}$
and $w \in {\Bbb C}$:
$$
\phi_0(z)=[1:z], \quad  \phi_1(w)=[w:1].
$$
The coordinates are connected by the relation $z^k w^l=1$.
Thus, $S^2(k,l)$ has two singular points $z=0$ and $w=0$
with the respective uniformizing groups ${\Bbb Z}_k$ and ${\Bbb Z}_l$.
Beyond the singular points, $S^2(k,l)$ possesses the structure
of a~complex manifold and, topologically, is a~two-dimensional sphere.

Now, consider the cotangent bundle $M_{k,l}=T^* S^2(k,l)$ and
study its structure.
In the tangent bundle $T ({\Bbb C}^2 \backslash \{0\})$,
we naturally distinguish the subbundle constituted by the vertical tangent
vectors with respect to the action of~${\Bbb C}^*$
(the vertical vectors are those tangent to the orbits of the action).
In the cotangent bundle
$T^*({\Bbb C}^2 \backslash \{0\})= {\Lambda^{1,0} {\Bbb C}^2} \times ({\Bbb C}^2 \backslash\{0\})$
consider the subbundle~$E$ constituted by the covectors vanishing on the vertical vectors.
It is easy to see that
$$
E=\{(z_0 (l z_2 d z_1 - k z_1 d z_2), z_1, z_2 )\mid  z_0 \in {\Bbb C},
(z_1, z_2) \in {\Bbb C}^2 \backslash \{0\} \}.
$$
Then the action of ${\Bbb C}^*$
on ${\Bbb C}^2 \backslash\{0\}$
induces the action of~${\Bbb C}^*$ on~$E$ which has the following structure:
$$
\lambda \in {\Bbb C}^*: (z_0, z_1, z_2) \mapsto (z_0 \lambda^{k+l}, z_1
\lambda^{-k}, z_2 \lambda^{-l}).
$$
Since the quotient space of~$E$ by the action of~${\Bbb C}^*$
coincides obviously with~$T^*S^2(k,l)$, we thereby find that
$M_{k,l}$ can be identified with ${\Bbb C} P^2_- (k+l,-k,-l)$.
In this event, the projective line~$S^2(k,l)$ is embedded into~$M_{k,l}$
as a~complex submanifold with singularities $\{ z_0=0 \}$.

Thus, $M_{k,l}$ is a~complex orbifold with two singular
points and the uniformizing groups ${\Bbb Z}_k$ and ${\Bbb Z}_l$ at these points.
In particular, if one of the numbers $k$ and $l$ equals ~1
then there is only one singularity, and if $k=l=1$ then we obtain
the cotangent bundle over the standard two-dimensional sphere without singularities.
Consider the two charts with local coordinates
$(z,\alpha) \in {\Bbb C}^2$ and $(w, \beta) \in {\Bbb C}^2$
which determine an~atlas on~$M_{k,l}$:
$$
\psi_1(\alpha, z)=[\alpha: 1:z],\quad  \psi_2(\beta,w)=[\beta : w : 1].
$$
The coordinate systems are connected by the relations $z^k w^l=1$
and
$\alpha z=\beta w$.
Moreover, the uniformizing groups~${\Bbb Z}_k$ and~${\Bbb Z}_l$
acting in each coordinate system are presented
in the group~$SU(2)$ acting standardly on~${\Bbb C}^2$:
$$
%\gather
{\Bbb Z}_k=\left\{ \left( \begin{array}{cc} \omega & 0 \\
0 & \omega^{-1} \end{array} \right) | \omega \in {\Bbb C},\
\omega^k=1 \right\},
\\
{\Bbb Z}_l=\left\{ \left( \begin{array}{cc} \omega & 0 \\
0 & \omega^{-1} \end{array} \right) | \omega \in {\Bbb C},\
\omega^l=1 \right\}.
%\endgather
$$
The projective line~$S^2(k,l)$ is embedded into~$M_{k,l}$ as a~complex
submanifold (with singularities) ${\{ \alpha\,{=}\beta{=}\,0\!\}}$.

Below we need the analytic structure of~$M_{k,l}$ at infinity.
Therefore, consider the isomorphism of complex manifolds
$$
%\gathered
\begin{array}{c}
\tau_{k,l}: ({\Bbb C}^2/{\Bbb Z}_{k+l}) \backslash\{0\}={\Bbb C}
P^2_+(k+l,-k,-l)\backslash \{[1:0:0]\} \rightarrow M_{k,l}\backslash
S^2(k,l),
\\
\tau_{k,l} [z_0:z_1:z_2]=[z_0:z_1:z_2].
%\endgathered
%                                \tag1
\end{array}\eqno{(1)}
$$
It is easy to see that $\tau_{k,l}$ ``glues in'' the projective line~$S^2(k,l)$
instead of the origin in~${\Bbb C}^2/{\Bbb Z}_{k+l}$.
Moreover, different orbits of the action of~${\Bbb C}^*$ on
${\Bbb C}^2 \backslash \{0\}$ ``intersecting'' at the point
$0 \in {\Bbb C}^2$
do not intersect any longer in~$M_{k,l}$;
that is, we have an~analog of the blow-up operation.

In~[1] the author constructed Ricci-flat metrics on~$M_{k,l}$
by means of special coordinates $\rho, \theta, \phi, \psi$.
The orbifold $M_{k,l}$ appeared as a~cylinder of the bundle
of the lens space~$L(-1,k+l)$ over~$S^2(k,l)$.
Moreover, $(\theta, \phi, \psi)$ were the coordinates in~$L(-1,k+l)$,
while $\rho$ was the coordinate along the element of the cylinder.
The metrics looked as follows:
$$
\begin{array}{c}
ds^2=(\mbox{ch} \rho-a \cos \theta) \left( d\rho^2+d\theta^2
\right)+
\\
\frac{\mbox{sh}^2 \rho}{\mbox{ch} \rho-a \cos \theta} \left( d\psi +
\cos \theta d\phi \right)^2 + \frac{\sin^2 \theta}{\mbox{ch} \rho-a
\cos \theta} \left( a d \psi +\mbox{ch} \rho d \phi \right)^2,
\end{array}\eqno{(2)}
$$
where $a=\frac{l-k}{l+k}$. For $k=l=1$ (i.e., for $a=0$) the
metric~(2) coincides with the Eguchi--Hanson metric of
cohomogeneity~1. For $a \neq 0$ we obtain a~new metric of
cohomogeneity ~2.

It is easy to note that the space $L(-1,k+l)$ is a~spherical
subbundle in $T^*S^2(k,l)=M_{k,l}$ over~$S^2(k,l)$.
This enables us to establish a~connection between the coordinates
$(\theta, \rho, \psi,\phi)$ and $(z,\alpha)$ or $(w,\beta)$.
Namely, we consider the following change of coordinates:
$$
z=\frac{\sin \frac{\theta}{2}}{\bigl(\cos
\frac{\theta}{2}\bigr)^{l/k}} e^{-i l (a \psi+\phi)}, \quad
\alpha=k l \mbox{sh \ }{\rho\over 2} \left(\cos {\theta\over
2}\right)^{1+\frac{l}{k}} e^{i l (\psi+\phi)},
$$
$$
w=\frac{\cos \frac{\theta}{2}}
{\bigl(\sin \frac{\theta}{2}\bigr)^{k/l}}
e^{i k (a \psi+\phi)},
\quad
 \beta= k l \mbox{sh \ }{\rho\over 2} \left(\sin {\theta\over 2}
\right)^{1+\frac{k}{l}} e^{i k (\psi-\phi)}.
$$
The metric $ds^2$ becomes smooth upon this change in each of
two charts (the exact expressions for the metric in the variables
$(z,\alpha)$ and $(w, \beta)$ are rather bulky and not given here).

Recall that a~Riemannian metric on a~two-dimensional complex
manifold (orbifold) is called a~{\it special K\"ahler metric\/}
if its holonomy group lies in the group~$SU(2)$ presented standardly in the tangent space.
One of the main results of the article is the following assertion:

\vskip0.2cm

{\bf Theorem 1.} {\it The space $M_{k,l}$ with metric {\rm(2)} is
a~special K\"ahler orbifold.}

\vskip0.2cm

{\bf Proof.} Here is a~K\"ahlerian form~$\omega$ which agrees with
metric~(2):
$$
\omega=\mbox{sh \ }\rho\, d\rho \wedge d\psi - a \sin \theta
\,d\theta \wedge d\psi + \mbox{sh \ }\rho \cos \theta \,d\rho \wedge
d\phi - \mbox{ch \ }\rho \sin \theta \,d\theta \wedge d\phi.
$$

Straightforward calculations show that this form is smooth
outside the singular points and is closed (and consequently parallel).
Hence, we can immediately conclude that $M_{k,l}$ is
a~K\"ahler manifold and the holonomy group lies in~$U(2)$.
One well-known result of differential geometry
(for example, see~[9]) states that
a~smooth Ricci-flat K\"ahler metric
on a~simply connected manifold is
a~special K\"ahler metric; i.e., its holonomy group lies in~$SU(2)$.
This result is of local character; therefore,
applying it to our case, we can conclude that parallel translation
along small loops away from the pair of singular points lies in~$SU(2)$.
Among the other things, the uniformizing group of the orbifold
has a~contribution into translation ``around'' a~singular point.
However, the elements of the uniformizing groups as well lie in~$SU(2)$;
therefore, the holonomy group~$M_{k,l}$ coincides with~$SU(2)$.

The theorem is proven.

\section[]{Applications to the Geometry of~$K3$-Surfaces}

{\bf Kummer's construction}.
Recall the construction of a~$K3$-surface by means of Kummer's construction.
Consider a~complex two-dimensional torus
$T^4={\Bbb C}^2/\Lambda$, where
$\Lambda=\{ (a+i b, c+i d) \mid  a,b,c,d \in {\Bbb Z}\}$ is
a~lattice in ${\Bbb C}^2$.
Define the involution $\sigma: T^4\rightarrow T^4$ as follows:
$$
\sigma: (z_1,z_2)+\Lambda \mapsto (-z_1, -z_2) +\Lambda.
$$
The involution $\sigma$ has $16$ fixed points; namely, the points
$(z_1, z_2)$ for
$z_i \in \bigl\{ 0, \frac{1}{2}, \frac{1}{2} i, \frac{1}{2}+\frac{1}{2} i \bigr\}$.
If we consider a~flat metric on~$T^4$ then $\sigma$ becomes
an~isometry and $X=T^4/\sigma$ represents a~special
K\"ahler orbifold
with $16$ singular points. In this event, a~neighborhood of each point looks like
${\Bbb C}^2/\{\pm 1 \}$.
Now, consider the complex surface~$Y$ obtained by blowing up each
of the~$16$ singular points; this is a~$K3$-surface.

The blow-up construction can be carried out as follows:
Let $u=(u_1,u_2)$ be the Euclidean coordinates in~${\Bbb C}^2$;
consider a~spherical neighborhood
$\Delta = \{|u|^2 \leq \varepsilon \}/{\Bbb Z}_2 \subset {\Bbb C}^2/{\Bbb Z}_2$
in which the blow-up will happen.
Now, consider $M_{1,1}=T^*S^2$ and the mapping
$\tau_{1,1} : {\Bbb C}^2 \backslash \{0\} \rightarrow M_{1,1} \backslash S^2$
constructed in~(1).
The mapping $\tau_{1,1}$ induces a~complex isomorphism between the open
submanifold~$\Delta \backslash \{0\}$ and an~open submanifold in~$M_{1,1} \backslash S^2$.
Deleting the point $0 \in \Delta$ and carrying out the identification
by means of the indicated isomorphism, we obtain a~complex manifold
without a~singular point.

\medskip
{\bf Page's construction}. It is well known [4] that $Y$ possesses
a~$58$-dimensional family~${\cal S}$ of metrics with
holonomy~$SU(2)$. In~[3] Page proposed a~geometric description of
the moduli space of the metrics with holonomy~$SU(2)$. We briefly
describe his approach. Considering all flat metrics on~$T^4$, we
obtain the family ${\cal S}_2$ of metrics of holonomy~$SU(2)$ on the
orbifold~$X$. Neighborhoods of $16$ singular points of~$X$ are
isometric to neighborhoods of the origin in~${\Bbb C}/ \{\pm 1\}$.
For each singular point we cut its neighborhood $(0, \varepsilon_1)
\times S^3/\{\pm 1\}$. Then the ``collar'' of the boundary
$(\varepsilon_1, \varepsilon_2)\times S^3/\{\pm1\}$ is ``almost''
isometric to the open spherical layer $M_{1,1}$ with a~metric
homothetic to the Eguchi--Hanson metric with the homothety
coefficient~$t$. Decreasing~$t$, we can make the metrics arbitrarily
close on the collar; moreover, both metrics have holonomy~$SU(2)$.
Now, using the analytic tools, we can show that if a~neighborhood is
sufficiently small then deformation of both metrics gives a~smooth
metric on~$Y$ with holonomy~$SU(2)$. A~rigorous justification of
this construction was given later in~[10,\,11,\,5].

This approach to the moduli space of a~$K3$-surface gives
a~geometrically clear explanation of the dimension~$58$. Indeed, the
space ${\cal S}_2$ is ten-dimensional. Now, in a~neighborhood of
each singular point, the metrics in~${\cal S}_2$ possess the group
of isometries~$SO(4)$ on~$(\varepsilon_1, \varepsilon_2) \times
S^3/\{\pm 1\}$. The subgroup~$U(2) \subset SO(4)$ leaves the
Eguchi--Hanson metric unchanged. This yields a~family of metrics of
dimension $\dim (SO(4)/U(2))=2$. If we use the parameter~$t$, we
obtain a~three-dimensional family of different metrics in
a~neighborhood of each singular point. Thus, taking all singular
points into account, we conclude that the dimension of the whole
family equals $10+16\cdot 3=58$. Actually, Page proposed a~geometric
description of the moduli space~${\cal S}$ of the metrics of
holonomy~$SU(2)$ on~$Y$ in a~neighborhood of the limit family~${\cal
S}_2$.

\medskip
{\boldmath\bf Resolution of singularities of type ${\Bbb C}^2/{\Bbb Z}_{k+l}$ by means of~$M_{k,l}$}.
We propose to consider~$M_{k,l}$ as a~space which enables us
to resolve singularities of higher order by reducing them
to singularities of less order.
More exactly, consider a~singular point of some complex manifold
such that the manifold in a~neighborhood of this point looks like
${\Bbb C}^2/{\Bbb Z}_{k+l}$.
Let $u=(u_1,u_2)$ be the Euclidean coordinates in~${\Bbb C}^2$
and let $\Delta=\{ |u|^2 \leq \varepsilon \}/{\Bbb Z}_{k+l}$ be
a~neighborhood of the manifold in which the blow-up happens.
Now, consider $M_{k,l}$ and the mapping
$\tau_{k,l}: ({\Bbb C}^2/{\Bbb Z}_{k+l}) \backslash \{0\} \rightarrow M_{k,l} \backslash S^2(k,l)$
defined in~(1).
This mapping induces a~complex isomorphism between the manifold
$\Delta \backslash \{0\}$ and an~open submanifold in~$M_{k,l} \backslash S^2(k,l)$.
Removing the singular point~$0 \in \Delta$ and
carrying out the identification by means of the indicated isomorphism,
we obtain a~complex manifold which has two singular points of types
${\Bbb C}^2/{\Bbb Z}_k$ and ${\Bbb C}^2/{\Bbb Z}_l$
instead of one singular point of the form~${\Bbb C}^2/{\Bbb Z}_{k+l}$.
Repeating this procedure, we can
now resolve all resulting singular points.

Pursuing the goal to generalize Page's construction, consider the following question:
What groups
${\Bbb Z}_p \subset SU(2)$ for $p>2$ can act on
$T^4={\Bbb C}^2/\Lambda$ by isometries?
It is clear that after an~appropriate choice of a~unitary basis
we can assume that
$$
\Z_p=\left\{ \left(
\begin{array}{cc} \omega^q & 0 \\
0 & \bar{\omega}^q
\end{array} \right) | q \in \Z
\right\},
$$
where $\omega=e^{\frac{2 \pi}{p} i}$ is the primitive root of
degree~$p$ of unity. Consider a~nonzero element $\lambda=(\lambda_1,
\lambda_2) \in \Lambda$ of the lattice. Let $\Pi=\{(\lambda_1 z,
\lambda_2 \bar{z} ) \mid  z \in {\Bbb C} \}$ be a~real
two-dimensional plane invariant under~${\Bbb Z}_p$. Since $\Lambda$
is invariant under the action of~${\Bbb Z}_p$, $\Pi \cap \Lambda$ is
a~lattice containing the vectors~$\lambda$, $\lambda \omega$,
and~$\lambda \omega^2$. Consequently, there is a~polynomial of the
second degree with integer coefficients whose root is the primitive
root~$\omega$. Hence, this is a~polynomial of division of the disk
which leaves only three possibilities: $p=3,4,6$. However, the cases
$p=4,6$ mean the presence of singular points in~$T^4/{\Bbb Z}_p$
whose neighborhoods are not modeled with the cone~${\Bbb C}^2/{\Bbb
Z}_p$, and we cannot resolve them using our construction. We are
left with the only case $p=3$ which will be considered in detail.

Suppose that the flat metric on ${\Bbb C}^2$ is determined by
the real part of the standard Hermitian product and
$\Lambda_0= \{ a e_1 +b e_2 \mid  a,b \in {\Bbb Z} \}$,
where $e_1=1$ and $e_2=e^{\frac{\pi}{3} i}$.
The arguments above demonstrate that the general lattice~$\Lambda$
invariant under the action of~${\Bbb Z}_3$
has the form
$$
\Lambda=\{(\lambda_1 z_1+\lambda_2 z_2, \mu_1 z_1 + \mu_2 z_2 ) \mid  z_1, z_2 \in
\Lambda_0 \} \cong \Lambda_0 \oplus \Lambda_0,
$$
where $\lambda_1$, $\lambda_2$, $\mu_1$, and~$\mu_2$ are complex parameters
which $\Lambda$ depends on, and
$\lambda_1 \mu_2 - \lambda_2 \mu_1 \neq 0$.
The action of the group ${\Bbb Z}_3$ on $T^4={\Bbb C}^2/\Lambda$
is generated by the transformation
$$
\gamma: (z_1, z_2) +\Lambda \mapsto (z_1 e^{\frac{2 \pi}{3} i},
z_2 e^{-\frac{2 \pi}{3} i}) +\Lambda.
$$
The quotient space $X=T^4/{\Bbb Z}_3$ is a~K\"ahler orbifold with
nine singular points which correspond to the fixed (with respect
to~$\gamma$) points in~$T^4$. These points have the form $(\lambda_1
z_1+\lambda_2 z_2, \mu_1 z_1+ \mu_2 z_2 )$, where $z_i \in \bigl\{
0, \frac{1}{3} e_1 + \frac{1}{3} e_2, \frac{2}{3} e_1 + \frac{2}{3}
e_2 \bigr\}$. We denote by $\{ s_1, \dots, s_9 \}$ the singular
points in~$X$. Let ${\cal S}_3$ be the moduli space of the flat
metrics on~$X$ with holonomy~$SU(2)$ corresponding to all possible
values of the parameters~$\lambda_i$ and~$\mu_i$. The action of the
group $U(2) \subset GL_{\Bbb C}(2)$ leaves the metric on~${\Bbb
C}^2$ unchanged and consequently induces the trivial action
on~${\cal S}_3$. Therefore, $\dim ({\cal S}_3)=\dim(GL_{\Bbb
C}(2)/U(2))=4$.

Let $X'$ be the complex surface obtained by successive
resolution of~$X$ at the singular points by means of~$M_{1,2}$ and~$M_{1,1}$.
Namely, suppose that $B_i \subset T^4/{\Bbb Z}_3$,
$i=1, \dots, 9$, are open balls of radius~$\varepsilon$
centered at the singular points of~$X$
and
$B_i' \subset B_i$ are the closed balls of radius $\varepsilon' < \varepsilon$
centered at the same points.
Thus, we obtain a~system of concentric neighborhoods of the singular
points in~$X$: $s_i \in B_i' \subset B_i$.
Choose a~sufficiently small~$\varepsilon$ such that
$B_i \cap B_j=\emptyset$ for $i \neq j$.
It is clear that $B_i \backslash B_i'$ is diffeomorphic to
$(\varepsilon', \varepsilon) \times S^3/{\Bbb Z}_3$.
Consider the space $M_{1,2}$ with the metric $t^2 ds^2$, where
$ds^2$ is the metric~(2) for $a=1/3$.
In the space $M_{1,2}$ consider the collar
$\tau_{1,2}((\varepsilon', \varepsilon) \times S^3/{\Bbb Z}_3)$,
where the mapping $\tau_{k,l}$ is defined in~(1).
Then the mapping $\tau_{1,2}$ defined on the collar tends to an~isometry
as $t\rightarrow 0$.
Consider the orbifold~$Y$ obtained by identification of
$X\backslash (\bigcup\nolimits_{i=1}^9 B_i')$
with nine copies of
$\tau_{1,2}((0, \varepsilon)\times S^3/{\Bbb Z}_3)$
by means of the mapping $\tau_{1,2}$ bounded on the collar
$(\varepsilon', \varepsilon) \times S^3/{\Bbb Z}_3$.
On the glued domains
the metric $t^2 ds^2$ is arbitrarily close to a~locally flat metric on~$X$ as
$t \rightarrow 0$.
The orbifold $Y$ has nine singularities $s_1', \dots, s_9'$
and looks locally like~${\Bbb C}^2/{\Bbb Z}_2$ in their neighborhoods.
Similarly, we consider a~sufficiently small
$\delta'<\delta <\varepsilon'$ and a~system of concentric
neighborhoods $s_i' \in C_i'\subset C_i$ in~$Y$ of radii~$\delta'$ and~$\delta$.
Using the mapping $\tau_{1,1}$, we identify
$Y\backslash (\bigcup\nolimits_{i=1}^9 C_i')$
with nine copies of
$\tau_{1,1}((0,\delta) \times S^3/{\Bbb Z}_2) \subset M_{1,1}$
over the collar $(\delta', \delta) \times S^3/{\Bbb Z}_2$
and obtain a~complex surface $X'$ without singularities.
Moreover, in the glued domains the metric on~$Y$ is arbitrarily
close to the metric $u^2 ds'{}^2$ as $u, \delta \rightarrow 0$,
where $ds'{}^2$ is the metric~(2) on~$M_{1,1}$ for $a=0$.
Let $d\tilde{s}^2$ be the metric on $Y$ obtained by smoothing the metrics
on~$X$, $M_{1,2}$, and $M_{1,1}$ in the domains of the described identification.

Recall that ${\cal S}$ is the moduli space of the metrics with the
holonomy group~$SU(2)$ on a~$K3$-surface, ${\cal S}_3$ is the moduli
space of the flat metrics on~$X$ with the holonomy group~$SU(2)$.

\vskip0.2cm

{\bf Theorem 2.} {\it The surface $X'$ is a~$K3$-surface and
consequently ${\cal S}_3$ is a~limit space for~${\cal S}$.
A~sufficiently small neighborhood of~${\cal S}_3$ in~${\cal S}$
consists of the $58$-dimensional family of metrics obtained by small
deformation of the family of metrics~$d\tilde{s}^2$ constructed as
described above as $\delta, t, u \rightarrow 0$.}

\vskip0.2cm

We give a~proof of this theorem in the next section.

\section[]{Connection with Multi-Instantons and the Proof of
Theorem~2}

The easiest way to justify the above construction rigorously
is to use the connection of the constructed metrics with multi-instantons.
Namely, multi-instantons are metrics with the holonomy group~$SU(2)$
of the following form~[7]:
$$
ds^2=\frac{1}{U} (d \tau + \mbox{\mathversion{bold}$\omega$} \cdot d
{\bf x} )^2+U d {\bf x} \cdot d {\bf x},\eqno{(3)}
$$
where ${\bold x} \in {\Bbb R}^3$, the variable $\tau$ is periodic,
$$
U=\sum_{i=1}^s \frac{1}{|{\bf x}-{\bf x}_i|},
$$
$$
\mbox{rot} \mbox{\mathversion{bold}$\omega$}= \mbox{grad} U,
$$
and ${\bold x}_i$ is a~set of~$s$ different points in the Euclidean
space ${\Bbb R}^3$. A~multi-instanton~(3) is a~smooth Riemannian
metric on some four-dimensional manifold.

Let $d\sigma$ be the area form of the level surface of ~$U$
in ~${\Bbb R}^3$.
It is easy to verify that the form
$$
\omega'_1=\frac{d U}{|\mbox{grad} U |} \wedge (d \tau +
\mbox{\mathversion{bold}$\omega$} \cdot d {\bf x})+ U d \sigma
\eqno{(4)}
$$
is a~closed K\"ahlerian form which agrees with metric~(3).

Consider the limit case when there are only two points
${\bold x}_1=(-1,0,0)$ and ${\bold x}_2=(1,0,0)$;
moreover, the first has multiplicity~$l$ and
the second has multiplicity~$k$. Take the coordinates $(r,\theta, \phi,\psi)$ as follows:
$$
x_1=\mbox{ch \ }\rho \cos \theta, \quad x_2=\mbox{sh \ }\rho \sin
\theta \cos \psi, \quad x_3=\mbox{sh \ }\rho \sin \theta \sin \psi,
\quad \tau=(l+k) \phi.
$$
Straightforward calculations show that with these coordinates
~(3) coincides with ~(2) to within multiplication by a~constant;
i.e., the constructed metric on~$T^* {\Bbb C} P^1(k,l)$
is a~limit case of a~multi-instanton.

Now, we can prove Theorem~2.
Our proof is similar to the arguments of~[5]; therefore,
as far as it is possible we will try
to keep the corresponding notations.
Let $T^4$ be the flat torus with the above-defined action of the group~${\Bbb Z}_3$.
Consider the torus $T^7=T^3 \times T^4$
with a~flat metric and extend the action of~${\Bbb Z}_3$ to~$T^7$
by making it trivial on~$T^3$.
Then the set~$S$ of singular points in~$T^7/{\Bbb Z}_3$ is
the disjoint union of nine tori~$T^3$; moreover, a~neighborhood~$T$
of the set~$S$ is isometric to the disjoint union of nine copies of~
$T^3 \times B_\zeta^4/{\Bbb Z}_3$, where $B_\zeta^4$
are open balls of radius~$\zeta$ in~${\Bbb R}^4={\Bbb C}^2$
for an~appropriate constant $\zeta>0$.

Now, choose an~arbitrary $\varepsilon>0$ and consider
a~multi-instanton $ds^2(t)$ on the four-dimensional manifold~$M_t$
given by the following three points:
${\bold x}_1=(-4 t^2/3,0,0)$,
${\bold x}_2=(2 t^2/3,t^2 \varepsilon,0)$,
and
${\bold x}_3=(2 t^2/3,-t^2 \varepsilon,0)$.
Denote by~$U_t$ the corresponding potential.
Let $\overline{U}({\bold x})=3/|{\bold x}|$
be the potential of the center of gravity with multiplicity~3.
It is easy to verify that
$$
|\nabla^i (U_t({\bold x})- \overline{U} ({\bold x}))| =O(t^4)\eqno{(5)}
$$
for $|{\bold x}|> \zeta^2/16$ and all $i\geq 0$
as $t \rightarrow 0$, where
$\nabla^i$
is the set of partial derivatives of order~$i$.
The metric~$d\bar{s}^2$ constructed for the potential~$\overline{U}$
is isometric to the flat metric on~${\Bbb C}^2/{\Bbb Z}_3$;
moreover, each domain
$\{ |{\bold x}| \leq r^2 \}$ is isometric to the ball~$B^4_r$.
It follows from~(5) that
$$
|\nabla^i (ds^2(t)-d\bar{s}^2)|= O(t^4)
                            \eqno{(6)}
$$
for $|{\bold x}|>\zeta^2/16$ as $t \rightarrow 0$.

Now, cut each of the nine neighborhoods
$B^4_\zeta/{\Bbb Z}_3$ of the singular points in
$X=T^4/{\Bbb Z}_3$ and, instead of each of them,
glue the domain in~$M_t$ defined by the condition $|{\bold x}| \leq \zeta^2$.
We obtain a~smooth four-dimensional manifold~$X'$.
In~$X'$ we define three domains~$A$, $B$, and~$C$ as follows:
the domain $A$ is the union of the neighborhoods given by the condition
$|{\bold x}|\leq \zeta^2/9$ in each of the nine glued copies of~$M_t$;
the domain $C$ is the complement in~$X'$ of the glued neighborhoods
and is isometric to~$X\backslash (\bigcup\nolimits_{i=1}^9 B^4_\zeta)$;
and finally $B=X' \backslash (A\cup C)$.

Consider the K\"ahlerian forms~$\omega'_1(t)$ and~$\overline{\omega}_1$
of the metrics $ds^2(t)$ and $d\bar{s}^2$ described in~(4).
Since the metric $ds^2(t)$ has the holonomy group~$SU(2)$,
we have a~parallel complex volume form~$\mu(t)$
on the K\"ahler manifold $(M_t,\omega'_1(t))$.
Put
$$
\omega'_2(t)=Re (\mu(t)), \ \omega'_3(t)= Im (\mu(t)).
$$
Then the triple of the parallel K\"ahlerian forms~$\omega'_1(t)$,
$\omega'_2(t)$, and~$\omega'_3(t)$ determines the hypercomplex
structure on~$M_t$ with the metric~$ds^2(t)$. Similarly, we
construct a~triple of parallel constant forms $\overline{\omega}_1$,
$\overline{\omega}_2$, and~$\overline{\omega}_3$ which determine the
flat hypercomplex structure on~${\Bbb C}^2/{\Bbb Z}_3$ with the
metric~$d\bar{s}^2$.

It follows from~(4) and~(5) that
$$
|\nabla^k (\omega'_i(t)-\overline{\omega}_i)|= O(t^4),\quad  i=1,2,3,
$$
for $|{\bold x}| > \zeta^2/16$ as $t \rightarrow 0$.

Consider the union $B$ of the spherical layers.
There exist $1$-forms $\eta'_i(t)$ and $\bar{\eta}_i$ such that $\omega'_i(t)=d\eta'_i(t)$
and
$\overline{\omega}_i=d \bar{\eta}_i$; moreover,
$$
|\nabla^k (\eta'_i (t)-\bar{\eta}_i)|=O(t^4) \mbox{\ as }t
\rightarrow 0.
$$
Consider a~real smooth increasing function $u(r)$ defined on the interval~$[0,\zeta]$
and possessing the following properties:
$$
0 \leq u(r) \leq 1; \quad u(r)=0 \mbox{\ for\ } 0 \leq r \leq
\zeta/3; \quad
 u(r)=1 \mbox{\ for\ } \zeta/2 \leq r \leq \zeta.
$$
Put $\eta_i(t)=u \bar{\eta}_i+(1-u) \eta'_i(t)$ and
$\omega_i(t)=d\eta_i(t)$, $i=1,2,3$.
We have thus constructed a~triple of closed $2$-forms on~$X'$
which coincide with the forms~$\overline{\omega}_i$ in the domain~$C$
and with the forms $\omega'_i(t)$ in the domain~$A$.
Moreover, it is easy to see that
$$
|\nabla^k (\omega_i(t)-\overline{\omega}_i)|=O(t^4)
                                \eqno{(7)}
$$
in the domain~$B \cup C$ as $t \rightarrow 0$.

Recall the definition of the $G_2$-structure.
Define a~$3$-form~$\phi_0$ on the space~${\Bbb R}^7$
with the standard Euclidean metric and orientation as follows:
$$
\phi_0=y_1 \wedge y_2 \wedge y_7+ y_1 \wedge y_3 \wedge y_6+ y_1
\wedge y_4 \wedge y_5+ y_2 \wedge y_3 \wedge y_5
$$
$$
- y_2 \wedge y_4 \wedge y_6+ y_3 \wedge y_4 \wedge y_7+ y_5 \wedge
y_6 \wedge y_7,
$$
where $y_1,\dots,y_7$ is the standard orthonormal positively
oriented basis for $({\Bbb R}^7)^*$. Moreover, the dual $4$-form
with respect to the Hodge operator looks like:
$$
*\phi_0=y_1 \wedge y_2 \wedge y_3 \wedge y_4+ y_1 \wedge y_2 \wedge
y_5 \wedge y_6- y_1 \wedge y_3 \wedge y_5 \wedge y_7+ y_1 \wedge y_4
\wedge y_6 \wedge y_7
$$
$$
+ y_2 \wedge y_3 \wedge y_6 \wedge y_7+ y_2 \wedge y_4 \wedge y_5
\wedge y_7+ y_3 \wedge y_4 \wedge y_5 \wedge y_6.
$$
The subgroup in~$GL_+({\Bbb R}^7)$, preserving the form $\phi_0$
or~$*\phi_0$, coincides with the group~$G_2$.

Now, if $M$ is a~seven-dimensional oriented manifold then let
$\Lambda_+^3 M$ and $\Lambda_+^4 M$ be the subbundles in~$\Lambda^3 T^*M$
and~$\Lambda^4 T^* M$ constituted by the forms
that have the above-mentioned form~$\phi_0$ and~$*\phi_0$
at each point $p \in M$ in an~appropriate oriented basis~$T_p^*M$.
It is easily seen that $\Lambda_+^3 M$ and $\Lambda_+^4 M$ are
open subbundles in~$\Lambda^3 T^*M$ and $\Lambda^4 T^* M$,
and we obtain the natural identification mapping
$\Theta: \Lambda_+^3 M \rightarrow \Lambda_+^4 M$
which takes each form looking locally like~$\phi_0$
into a~form looking locally like~$*\phi_0$.
We say that a~section $\phi$ of the bundle $\Lambda_+^3 M$
{\it determines a~$G_2$-structure on~$M$.} This form~$\phi$
determines uniquely the Riemannian metric with respect to which
the identification operator~$\Theta$ becomes the Hodge operator~$*$.
Moreover, if the forms~$\phi$ and $*\phi$ are closed then
the $G_2$-structure is torsion-free and the holonomy group
of the Riemannian manifold $M$ is contained in~$G_2 \subset SO(7)$.

Defined a~flat $G_2$-structure~$\bar{\phi}$ on $T^3 \times X=T^7/{\Bbb Z}_3$
and its dual $* \bar{\phi}$ as follows:
$$
\bar{\phi}=\bar{\omega}_1 \wedge \delta_1+\bar{\omega}_2 \wedge
\delta_2+\bar{\omega}_3 \wedge \delta_3+\delta_1 \wedge \delta_2
\wedge \delta_3,
$$
$$
* \bar{\phi}=\bar{\omega}_1 \wedge \delta_2 \wedge
\delta_3+\bar{\omega}_2 \wedge \delta_3 \wedge
\delta_1+\bar{\omega}_3 \wedge \delta_1 \wedge \delta_2+\frac{1}{2}
\bar{\omega}_1 \wedge \bar{\omega}_1,
$$
where $\delta_1$, $\delta_2$, and~$\delta_3$ are constant
orthonormal $1$-forms on~$T^3$ extended to the whole
$T^7/{\Bbb Z}_3=T^3 \times X$.

\goodbreak

\noindent
Define the following $3$- and $4$-forms on~$M=T^3 \times X'$:
$$
\phi_t=\omega_1(t) \wedge \delta_1+\omega_2(t) \wedge
\delta_2+\omega_3(t) \wedge \delta_3+\delta_1 \wedge \delta_2
\wedge \delta_3,
$$
$$
v_t=\omega_1(t) \wedge \delta_2 \wedge \delta_3+\omega_2(t)
\wedge \delta_3 \wedge \delta_1+\omega_3(t) \wedge \delta_1
\wedge \delta_2+\frac{1}{2} \omega_1(t) \wedge \omega_1(t).
$$
It is clear that these forms coincide with~$\bar{\phi}$ and $*\bar{\phi}$
respectively in the domain~$T^3 \times C$.
Since all forms~$\omega_i(t)$ and $\delta_i$ are closed,
the forms~$\phi_t$ and $v_t$ are closed on~$M$ too.

In the domains $T^3 \times A$ and $T^3 \times C$, the triple of the forms~$\omega_i(t)$
is a~triple of K\"ahlerian forms determining a~hypercomplex structure;
therefore, the form $\phi_t$ determines a~torsion-free $G_2$-structure
and $v_t=\Theta(\phi_t)$.
In the domain~$T^3 \times B$, the triple~$\omega_i(t)$
does not determine a~hypercomplex structure in general,
thereby we cannot even guarantee a~priori that $\phi_t$ determines
a~$G_2$-structure.
However, $\Lambda_+^3 (M)$ is open in $\Lambda^3 T^* (M)$;
therefore, it follows from~(7) that for a~sufficiently small~$t$
we have $\phi_t \in C^\infty(\Lambda_+^3 (M))$.
Moreover, the form~$v_t$ differs from~$\Theta(\phi_t)$ in general.
Define the $3$-form $\psi_t$ on~$M$ by the relation
$*\psi_t=\Theta(\phi_t)-v_t$, where the Hodge operator is
defined with respect to the Riemannian metric $g$ given
by the $G_2$-structure~$\phi_t$.
It is obvious that $d^* \psi_t=d^*\phi_t$.

The following theorems are proven in~[5]:

\vskip0.2cm

{\bf Theorem A.} {\it Let $E_1,\dots, E_5$ be positive constants.
Then there are positive constants~$\kappa$ and $K$ depending on
$E_1, \dots, E_5$ and such that the following property holds for
every $0 < t <\kappa$.

Let $M$ be a~compact seven-dimensional manifold and let $\phi$ be a~smooth closed
form in~$C^{\infty}\bigl(\Lambda_+^3M\bigr)$.
Suppose that $\psi$ is a~smooth $3$-form on~$M$ such that
$d^*\psi=d^*\phi$  and the following are fulfilled:

{(i)} $\| \psi\|_2 \leq E_1 t^4$ and $\|\psi\|_{C^{1,1/2}}\leq E_1
t^4$;

{(ii)} if $\chi \in C^{1,1/2}(\Lambda^3 T^*M)$ and $d \chi=0$ then
$$
 \|\chi\|_{C^0}\leq E_2(t \|\nabla\chi\|_{C^0}+t^{-7/2}\|\chi\|_2),
$$
$$
\|\nabla \chi\|_{C^0}+t^{1/2}[\nabla \chi]_{1/2}\leq
E_3(\|d^*\chi\|_{C^0}
+ t^{1/2}[d^*\chi]_{1/2}+t^{-9/2}\|\chi\|_2);
$$

{(iii)} $1 \leq E_4 \mbox{vol} (M)$;

{(iv)} if $f$ is a~smooth real function and $\int\nolimits_M f
d\mu=0$ then $\|f\|_2\leq E_5 \|df\|_2$.

Then there is $\eta \in C^\infty(\Lambda^2 T^* M)$ such that
$\|d\eta\|_{C^0} \leq K t^{1/2}$ and $\tilde{\phi}=\phi+d\eta$ is
a~smooth torsion-free $G_2$-structure on~$M$.}

\vskip0.2cm

{\bf Theorem B.} {\it Let $D_1,\dots, D_5$ be positive constants.
Then there are positive constants $E_1,\dots,E_5$ and $\lambda$
depending only on $D_1,\dots, D_5$ such that the following property
holds for each $t \in (0,\lambda]$.

Let $M$ be a~compact seven-dimensional
manifold and let $\phi$ be
a~closed form in~$C^\infty\bigl(\Lambda_+^3M\bigr)$.
Let $g$ be the metric associated with~$\phi$.
Suppose that $\psi$ is a~smooth $3$-form on~$M$ such that
$d^*\phi=d^*\psi$ and the following are fulfilled:

{(i)} $\|\psi\|_2 \leq D_1 t^4$ and $\|\psi\|_{C^2}\leq D_1 t^4$;

{(ii)} the injectivity radius~$\delta(g)$ satisfies the inequality
$\delta(g) \geq D_2 t$;

{(iii)} the Riemannian tensor $ R(g)$ of the metric $g$ satisfies
the inequality $\|R(g)\|_{C^0}\leq D_3 t^{-2}$;

{(iv)} the volume $\mbox{vol}(M)$ satisfies the inequality
$\mbox{vol}(M) \geq D_4$;

{(v)} the diameter~$\mbox{diam} (M)$ satisfies the inequality
$\mbox{diam}(M) \leq D_5$.

Then the conditions~{\rm(i)--(iv)} of Theorem~{\rm A} are satisfied
for~$(M,\phi)$.}

\vskip0.2cm

%\goodbreak

We want to apply these theorems to the form $\phi_t$ on~$M$ with
the associated metric~$g$.
Indeed, (7) implies validity of condition~(i),
while condition (iv) follows trivially from the construction.
Now, note that $ds^2(t)=t^2\,ds^2(1)$.
Hence, we find that the injectivity radius grows linearly with
the increase of~$t$ which proves~(ii) and~(iii).
The same arguments imply boundedness of the diameter, i.e.,~(v).
Hence, there is a~torsion-free $G_2$-structure~$\tilde{\phi}$ close to~$\phi_t$.
Let $g'$ be the associated metric with the holonomy group~$G_2$ on~$M=T^3 \times X'$.

Since $X$ is simply connected, $X'$ is simply connected too.
Therefore, $\pi_1 (T^3 \times X')={\Bbb Z}^3$
and from~[5] we can conclude that the holonomy group~$(M,g')$
is equal to~$SU(2) \subset G_2$.
Hence, the metric on~$M$ is the direct product of a~flat metric on~$T^3$
and the metric $ds'{}^2$ with the holonomy group~$SU(2)$ on~$X'$.
In particular, $X'$ is a~$K3$-surface.

Now, we find out how the metric $ds'{}^2$ looks like near the remote singular points,
i.e., in the domain~$A$.
The metric~$g$ on~$T^3 \times X'$ close to~$g'$ in the domain~$A$
is the direct product of the flat metric on~$T^3$ and
the multi-instanton $ds^2=ds^2(t)$.
The multi-instanton $ds^2$ tends to a~flat metric on~${\Bbb C}^2/{\Bbb Z}_3$
as $t \rightarrow 0$ and to the metric~(2) on~$M_{1,2}$ as $\varepsilon \rightarrow 0$.
Hence, as $\varepsilon \rightarrow 0$,
the metric $ds^2$ is obtained from~$X$ by the resolution of the singular points
by means of~$M_{1,2}$ described in Theorem~2.
Now, if we consider a~small neighborhood of the points
${\bold x}_2$ and~${\bold x}_3$ and choose $\varepsilon$ so small that
the contribution of~${\bold x}_1$ in the potential $U_t$ is small
as compared with the contribution of~${\bold x}_2$ and~${\bold x}_3$
then the metric~$ds^2$ in this neighborhood is close to the metric~(2) on~$M_{1,1}$.
Thus, the metric $ds^2$ is obtained by the double resolution of singular points in~$X$
indicated in Theorem~2, while the metric~$ds'{}^2$ is its small deformation.

Estimate the dimension of the family of metrics constructed above.
In the process of resolution of the singularities~$s_i$ the freedom
of gluing of~$M_{1,2}$ is determined by the group~$U(2)$ which does
not change the complex structure on~$T^4$. However, the metric
on~$M_{1,2}$ has the group of isometries~$U(1)\times U(1) \subset
U(2)$; therefore, if we use the presence of the parameter~$t$
responsible for homothety then we obtain a~family of different
metrics with holonomy~$SU(2)$ of dimension~$3$ in a~neighborhood of
each point~$s_i$. In the processes of resolution of the
singularities~$s_i'$, as in Page's method, we also obtain a~family
of dimension~$3$. The dimension of~${\cal S}_3$ is equal to~$4$;
therefore, summing up the dimensions, we conclude that the dimension
of~${\cal S}$ in a~neighborhood of~${\cal S}_3$ equals~$58$; this is
exactly the dimension of the moduli space of the metrics with
holonomy~$SU(2)$ on a~$K3$-surface.

The proof of the theorem is complete.


\begin{thebibliography}{MMM}


\bibitem{no1}
Baza\u \i kin~Ya.~V. On some Ricci-flat metrics of cohomogeneity two
on complex line bundles. Siberian Math.~J. 2004. 45, 3, 410--415

\bibitem{Eguchi-Hanson}
T.~Eguchi, A.~J.~Hanson. Asymptotically flat self-dual solutions to
euclidean gravity. Physics Letters, {\bf 74B} (1978), No. 4,
249--251.

\bibitem{Page}
D.~N.~Page. A physical picture of the $K3$ gravitational instanton.
Physics Letters, {\bf 80B} (1978), No. 1,2. 55--57.

\bibitem{Yau}
S.-T.~Yau. On the Ricci curvature of a compact K\"ahler manifold and
the complex Monge-Amp\`{e}re equations. I. Communications on pure
and applied mathematics. {\bf 31} (1978). 339--411.

\bibitem{Joyce1}
D.~D.~Joyce. Compact Riemannian $7$-manifolds with holonomy $G_2$. I
and II. J. Differentional Geometry. {\bf 43} (1996), No. 2. 291--328
and 329--375.

\bibitem{Joyce2}
D.~D.~Joyce. Compact $8$-manifolds with holonomy $Spin(7)$. Inv.
Math. {\bf 123} (1996). 507--552.

\bibitem{Gibbons-Hawking}
G.~W.~Gibbons, S.~W.~Hawking. Gravitational multi-instantons.
Physics letters, {\bf 78B}, (1978), No. 4, 430--432.

\bibitem{Hitchin}
N.~J.~Hitchin. Polygons and gravitons. Math. Proc. Camb. Phil. Soc.
{\bf 85} (1979), 465--476.

\bibitem{Besse}
A.~Besse. Einstein manifolds. Springer-Verlag, NY, 1987.

\bibitem{Topiwala}
P.~Topiwala. A new proof of the existence of K\"ahler-Einstein
metrics on K3. I. Inventiones mathematicae. {\bf 89} (1987).
425--448.

\bibitem{Lebrun-Singer}
C.~LeBrun, M.~Singer. A Kummer-type construction of self dual
$4$-manifolds. Mathematische Annalen. {\bf 300} (1994). 165--180.

\end{thebibliography}
\end{document}